\documentclass{article}
\usepackage{geometry}
\usepackage[utf8]{inputenc}
\usepackage{amsmath,amsfonts,amssymb,amsthm,tabto,commath,multicol}
\theoremstyle{theorem}
\newtheorem{lemme}[section]{Lemma}
\theoremstyle{theorem}
\newtheorem{thm}[section]{Theorem}
\theoremstyle{theorem}
\newtheorem{prop}[section]{Proposition}
\newcommand*{\onto}{\ensuremath{\twoheadrightarrow}}
\newcommand{\xto}[1]{\xrightarrow{#1}}
\begin{document}

\begin{center}
{\Large \textbf{Quasi-automatic groups are asynchronously automatic}}

\begin{large}
\textsc{Benjamin Blanchette}\\ 
\end{large}
\end{center}

Classes of groups and semigroups have been defined in regards to their computation complexity in the past, notably by Epstein et al. with their early work on automatic groups, see [1]. More general classes have followed. First, the asynchronous case for groups was defined and studied, notably by the same group of authors. Afterwards, monoids and semigroups have been investigated, notably by Campbell et al., see [4], and Duncan et al., see [5].  However, a key feature of automatic groups does not hold for automatic semigroups; while the notion of automaticity for finitely generated groups does not depend on the choice of generators, it does for finitely generated semigroups. Showing there exists no automatic structure for a given generating set is doable, but doing so for any generating set can be tedious, see for example [6], where Hoffmann and Thomas build a finitely generated commutative semigroup which is not automatic under any choice of generating set. This makes the feature extremely valuable. Recently, Blanchette et al. introduced quasi-automatic semigroups, a more general class of semigroups for which the feature actually holds, see [2].\\

Several relations between these classes have been established. For instance, we know that the class of automatic groups sits strictly inside the class of asynchronously automatic groups, which itself is contained into the class of quasi-automatic groups. Geometric characterizations have also been made. Automatic groups are exactly those having the so-called Lipschitz property. Asynchronously automatic groups are characterized by the Lipschitz Hausdorff property of Epstein et al., along with a special function named a departure function. For quasi-automatic groups, we have a different weak form of the Lipschitz property that again characterizes the class geometrically.\\

Our main result here is to establish a new link between those classes; we will show that for groups, being quasi-automatic is equivalent to being asynchronously automatic. This is particularly useful because quasi-automatic structures are somewhat natural, as opposed to the more complex asynchronously automatic structures.\\

\begin{center}
{\Large \textbf{Context and definitions}}
\end{center}

Let $A$ be a generating set of a semigroup $S$ and $p:A^+ \onto S$ be the canonical map which maps a word to the element of $S$ it represents. We write $u\sim v$ if $u$ and $v$ represent the same element of $S$. A rational language $L \subseteq A^+$ that maps onto $S$ is called a \textit{dictionary}. A dictionary is called \textit{quasi-automatic }if for each $a$ in $A\cup \{\varepsilon\}$, the relation ${R_a^L=\{(u,v)\in L\times L| ua\sim v\}}$ is rational. If these relations are also recognizable by two-tape automata, this dictionary is called \textit{asynchronously automatic}. We use a theorem of Nivat (see [3], proposition 4) to work with these rational relations. The theorem states that a relation $R$ is rational if and only if there exist some finite $B$, some rational langage $H\subseteq B^*$ and some alphabetic morphisms $\alpha,\beta:B^*\to A^*$  verifying $(\alpha\times\beta)(H)=R$ such that for all $b\in B$, $\alpha(b)= \varepsilon$ or $\beta(b)=\varepsilon$, but not both. Such a pair of morphisms is called a \textit{Nivat bimorphism}.\\

A dictionary for a group is said to be \textit{Lipschitz Hausdorff }if whenever $u,v\in L$ are at most at distance $1$ in the Cayley graph of $G$, the Hausdorff distance between their respective sets of prefixes is uniformly bounded.  This can be rephrased as follows: there exists some $k\in \mathbb{N}$ such that for all $u,v\in L$ with $d(u,v)\leq 1$ and any prefix $p$ of $u$, there exists some prefix $q$ of $v$ such that $d(p,q)\leq k$. A dictionary for a group is said to be \textit{weakly Lipschitz }if whenever $u,v$ are at most at distance $1$ in the Cayley graph, we can rewrite $u=a_1…a_n$, $v=b_1…b_n$ where $a_i,b_i\in A \cup \{\varepsilon\}$ in such way that $d(a_1…a_i,b_1…b_i)$ is bounded by some $k\in \mathbb{N}$ for all $i\leq n$. A \textit{departure function }for a dictionary $L$ of a group is a function $D: \mathbb{N}\to \mathbb{N}$ such that if $xyz \in L$ and $|y|>D(n)$, then the distance between the group element $y$ represents and the neutral element of $G$ is at least $n$. This value is called the \textit{norm} of $y$ relative to $G$ and is written $|y|_G$.\\

Note that a weakly Lipschitz dictionary is also Lipschitz Hausdorff; given a prefix $p$ or $u$, there is an $i$ such that $a_1…a_i=p$, which is at most at distance $k$ from $b_1…b_i$, which is a prefix of $v$. \\

\begin{center}
{\Large \textbf{Main result}}
\end{center}

\begin{thm}
A finitely generated group $G$ is quasi-automatic if and only it is asynchronously automatic.
\end{thm}
For semigroups in general, asynchronously automatic dictionaries are quasi-automatic (see [2], Proposition 3.4.2). This implies the if part.  For the only if part, we use a theorem due to Epstein and Levy  (see [1], theorem 7.2.8) which goes as follows: a group is asynchronously automatic if it has a Lipschitz Hausdorff dictionary along with a departure function for that dictionary. This means it suffices to show there exist such a dictionary for any quasi-automatic group. This will be Proposition \textbf{4}; we first need a couple of lemmas.\\

\begin{lemme}
Let $L\subseteq A^+$ be a quasi-automatic dictionary of some semigroup $S$. There exist $k\in \mathbb{N}$, a rational dictionary $K\subseteq L$, some finite set $B$, some rational language $H\subseteq B^+$ with a Nivat bimorphism $\alpha\times \beta$ such that $(\alpha\times\beta)(H)=R_\varepsilon^K$ and for all $p,f,s\in B^*$,
$$pfs\in H,|f|>k \implies \alpha(f)\neq \varepsilon\neq \beta(f)$$
\end{lemme}
\textit{Proof} By definition, $R_\varepsilon^L$ is rational and by Nivat's theorem there exist some finite $B$, rational langage $H_1\subseteq B^*$ and alphabetic morphisms $B^*\to A^*$ such that $(\alpha\times\beta)(H_1)=R_\varepsilon^L$. Let $W=(Q,T,q_0,F)$ be a finite state deterministic automaton which recognizes $H_1$ and let $k$ be the number of states of $W$.  For each state $q\in Q$, we consider the set of circuits $y$ from $q$ to $q$ such that $\alpha(y)=\varepsilon$ and $0<|y|\leq k$. Note that this set finite. Let $L(q)$ be the set of accepted words of $W$ with paths including such a circuit. Explicitly, we have
$$L(q)=\{x\in B^+| q_0x=q\} \{y\in B^+| qy=q,\alpha(y)=\varepsilon,|y|\leq k\}\{ z\in B^+|qz\in F\}.$$
It is rational by concatenation of rational langages. Symmetrically, we define  $R(q)$ to be the set of accepted words whose path contains a non-empty circuit $y$ of length at most $k$ such that $\beta(y)=\varepsilon$. It is rational by the same reasoning.\\

We claim $K=L-\bigcup\limits_{q\in Q} \beta(L(q))\cup \alpha(R(q))$ is a rational dictionary of $S$. The relations ${R_a^K=R_a^L\cap (K\times K)}$ are rational by intersection. We show that $p(K)=S$.  Consider some $s\in S$ and $w\in L$ a word of minimal length such that $p(w)=s$.  Suppose for a contradiction that $w \notin K$. Then either $w\in \beta(L(q))$ or $w\in \alpha(R(q))$. Suppose the former; the other case is identical. Then $w=\beta(xyz)$ with $y$ a short circuit around $q$ with $\alpha(y)=\varepsilon$, $q_0x=q$ and $qz\in F$. Then $q_0xz=qz\in F$, indicating $xz\in H_1$. This implies $(\alpha\times\beta)(xz)\in R_\varepsilon^L$, which implies $\alpha(xz)\sim\beta(xz)$. Similarly, $\alpha(xyz)\sim \beta(xyz)$. Then
$$\beta(xyz)\sim \alpha(xyz)=\alpha(x)\alpha(y)\alpha(z)=\alpha(x)\varepsilon \;\alpha(z)=\alpha(xz)\sim \beta(xz)$$
We conclude that $\beta(xz)$ and $\beta(xyz)$ represent the same element of $S$. Since $\alpha\times \beta$ is a Nivat bimorphism, $\alpha(y)=\varepsilon$ implies $\beta(y)\neq \varepsilon$. This means $\beta(xz)$ is shorter than $w=\beta(xyz)$, a contradiction to the minimality of $w$.\\

Now let $H=(\alpha\times\beta)^{-1}(K\times K) \cap H_1$. Reverse morphisms preserve recognisability in any monoid; $H$ is therefore recognisable. But as a subset of the free monoid $B^*$, this means it is rational.  Furthermore, we have 
$$(\alpha\times\beta) (H)=(\alpha\times \beta)((\alpha\times\beta)^{-1}(K\times K)\cap H_1)=(K\times K) \cap (\alpha\times \beta)(H_1)= (K\times K) \cap R_\varepsilon^L=R_\varepsilon^K$$

We now show $H$ satisfies the lemma. Let $w=pfs\in H$ with $|f|>k$. Suppose for a contradiction that $\alpha(f)=\varepsilon$; here again, the symmetric case is identical. Since $H\subseteq H_1$, $pfs$ is recognised by $W$. Since the length of $f$ is greater than the number of states of $W$, its path during the reading of $w$ must contain some circuit $y$ of length at most $k$ around some $q$. Also, as $y$ is a factor of $f$, we know $\alpha(y)=\varepsilon$. As $w$ is an accepted word with path including an appropriate circuit $y$, $w\in L(q)$ and in turn $\beta(w)\not\in K$. This is a contradiction because $w\in H$ implies $ (\alpha\times \beta)(w)\in R_\varepsilon^K$ and thus $\beta(w)\in K$. \\
\newpage
\begin{lemme}
With the notation of the previous lemma, there exists $\ell \in \mathbb{N}$ such that 
$$(u,v)\in R_\varepsilon^K\implies |u|\leq \ell|v|.$$
\end{lemme}
\textit{Proof} Take $w\in H$ such that $(\alpha\times\beta)(w)=(u,v)$. Decompose $w$ into factors of length $k$. This gives us $w=f_1...f_{h+1}$ with $|f_i|=k$ for all $i\leq h$ and $|f_{h+1}|\leq k$. By the previous lemma, $\alpha(f_i)\neq \varepsilon\neq \beta(f_i)$ for all $i\leq h$.  We have $|u|=|\alpha(f_1...f_{h+1})|=|\alpha(f_1)|+...+|\alpha(f_{h+1})|$ and thus $h \leq |u| \leq kh+k$. Symmetrically, we have $h\leq |v|\leq kh+k$. Then $|u|\leq kh+k \leq k|v|+k\leq 2k|v|$ and $\ell=2k$ satisfies the condition of the lemma.\\

\begin{prop}
Let $L\subseteq A^+$ be a quasi-automatic dictionary for a group $G$. There exists a Lipschitz Hausdorff dictionary $K \subseteq L$ with a departure function.
\end{prop}
\textit{Proof} Start by taking $K\subseteq L$ like in the previous lemma and let $W=(Q,T,q_0,F)$ be a new finite state automaton recognising $K$; assume that $W$ is accessible and coaccessible. Let $c$ be the number of states of $W$ and $\ell$ be the constant given by the previous lemma. For all triplets $(q,q',g)$ with $q,q'\in Q$ and $g\in G$, set ${m(q,q',g)=\text{min}\{|w|:qw=q',p(w)=g\}}$
if there exist some $w$ with $qw=q'$ and $p(w)=g$, and $m(q,q',g)=0$ otherwise. We define $D:\mathbb{N}\to \mathbb{N}$ by setting 
$$D(n)=2c\ell +\ell \text{ max}\{m(q,q',g) |  q,q'\in Q,|g|_G\leq n\}$$

We claim this is a departure function for $K$.  Choose $xyz\in K$ such that $|y|>D(n)$ and suppose for a contradiction that $|y|_G\leq n$. Consider the path $q_0\xto{x}q_1\xto{y}q_2\xto{z}q_3\in F$. Choose a word $y'$ such that $q_1y'=q_2$ and $y'\sim y$ of minimal length $m(q_1,q_2,p(y))$. Since $W$ is accessible, coaccessible and has $c$ states, we can choose $x'$ and $z'\in A^*$ such that $|x'|,|z'|\leq c$ and $q_0x'=q_1$ and $q_2z'\in F$. This implies $x'yz'$ and $x'y'z'$ are in $K$. Also, since $y\sim y'$, we have $x'yz'\sim x'y'z'$. These are sufficient conditions to know that $(x'yz',x'y'z')\in R_{\varepsilon}^K$. The previous lemma then tells us that $|x'yz'|\leq \ell |x'y'z'|$.  Then,
$$|y|\leq |x'yz'|\leq \ell |x'y'z'| = \ell (|x'|+|y'|+|z'|)\leq \ell (c+|y'|+c)=2c\ell +\ell |y'|\leq D(n)$$

which is a contradiction to the choice of $y$. \\

\textit{Proof of Theorem 1} Let $G$ be a quasi-automatic group and $K$ be the quasi-automatic dictionary with departure function built in Lemma \textbf{2}.  By a theorem of Blanchette et al. (see [2], Theorem 4.5), as with all quasi-automatic dictionaries, it is weak Lipschitz. It is therefore Lipschitz Hausdorff by the remark at the end of the introduction. Applying the theorem of Epstein and Levy let us conclude that $G$ is asynchronously automatic.\\
\newpage
\begin{center}
{\Large \textbf{Conclusion}}
\end{center}
While we conjecture that quasi-automatic semigroups strictly contains the class of asynchronously automatic semigroups, we don't have an explicit example of a quasi-automatic semigroup which is not asynchronously automatic. On the level of relations, the inclusion is known to be strict. Rational relations are stable under union while relations recognised by two-tape automata are not; this indicates how different they are.  Building examples that are asynchronously automatic but not automatic is far from trivial (see [1] example 7.4.1, the Baumslag--Solitar group). As a similar problem, it hints us that finding an appropriate quasi-automatic semigroup and showing it is not asynchronously automatic could be a hard task.\\

\begin{center}
\section*{References}
\end{center}

\begin{small}
[1] D. Epstein et al., \textit{Word processing in groups}, Jones and Bertlett Publishers, Boston, 1992.

[2] B. Blanchette, C. Choffrut and C. Reutenauer, \textit{Quasi-automatic semigroups}, Theoretical Computer Science, 2019. In Press. DOI:  10.1016/j.tcs.2019.01.002

[3] M. Nivat, \textit{Transduction des langages de Chomsky}, Annales de l'Institut Fourier 18, 339-456, 1968.

[4] C. Campbell et al., \textit{Automatic semigroups}, Theoretical Computer Science 250, 365-391, 2001.

[5] A. Duncan, E. Robertson, N. Ruskuc, \textit{Automatic monoids and change of generators}, Mathematical Proceedings of the Cambridge Philosophical Society 127, 403-409, 1999.

[6] M. Hoffman, R.M. Thomas, \textit{Automaticity and commutative semigroups}, Glasgow Mathematical Journal 44, 167-176, 2002.

\end{small}
\end{document}